\documentclass[12pt]{article}
\usepackage[all]{xy}
\usepackage{amscd, amssymb,latexsym,amsmath, amscd, amsmath}

\numberwithin{equation}{section}
\newtheorem{theorem}{Theorem}[section]

\newtheorem{lemma}[theorem]{Lemma}

\newtheorem{remark}[theorem]{Remark}
\def\Tor{\text{\rm Tor}}
\def\Hom{\text{\rm Hom}}
\def\ubeta{\boldsymbol\beta}
\def\ualpha{\boldsymbol\alpha}

\newcommand{\C}{{\mathbb C}}
\newcommand{\Q}{{\mathbb Q}}

\newcommand{\Z}{{\mathbb Z}}

\begin{document}

\begin{center}
{\Large \bf Descent of line bundles to GIT quotients of flag varieties by maximal torus} \\[6mm]
{\bf Shrawan Kumar}
\end{center}
\vskip 1cm

\section{Introduction}
Let $G$ be a connected, simply-connected  semisimple complex algebraic group with a maximal torus
$T$ and let $P$ be a parabolic subgroup containing $T$. We denote their Lie algebras by the corresponding Gothic characters.
The following theorem is one of our main results.

\paragraph*{Theorem 3.9} {\it Let ${\cal L}_P (\lambda)$ be a homogeneous ample line bundle on the flag variety $Y=G/P$.  Then, the line bundle ${\cal L}_P (\lambda)$
descends to a line bundle on the GIT quotient $Y (\lambda) // T $  if
and only if for all the semisimple subalgebras $\mathfrak{s}$ of $\mathfrak{g}$ containing $\mathfrak{t}$ (in particular,  rank $\mathfrak{s}$ = rank
 $\mathfrak{g})$,
$$\lambda \in \displaystyle{\sum_{\alpha \in \triangle^+ (\mathfrak{s})}} \Z \alpha,$$
where $\triangle^+ (\mathfrak{s})$ is the set of positive roots of $\mathfrak{s}$.}

Using the above theorem, we explicitly get exactly which 
line bundles  ${\cal L}_P (\lambda)$ descend to 
 $Y(\lambda) // T$. This is our second main result.

In the following $Q$ (resp., $\Lambda$) is the root (resp., weight) lattice and we follow the indexing
convention as in Bourbaki [B].

\paragraph*{Theorem 3.10}  {\it With the notation as in the above theorem, the
line bundle ${\cal L}_P (\lambda)$ descends to a line bundle on the GIT quotient
$Y (\lambda) //T$ if and only if $\lambda$ is of the following form depending upon the type of $G$.}

\begin{enumerate}
\item[a)] $G$ of type $A_{\ell}\, (\ell \geq 1): \lambda \in Q$.
\item[b)]$G $ of type $B_\ell \,(\ell \geq 3): \lambda \in 2Q$.
\item[c)]$G $ of type $C_\ell \,(\ell \geq 2): \lambda \in \Z 2\alpha_1+
\cdots +\Z 2 \alpha_{\ell-1}+ \Z \alpha_\ell = 2\Lambda $.
\item[$d_1$)]  $G $ of type $D_{4}: \lambda \in \{n_1 \alpha_1+ 2n_2 \alpha_2+
n_3 \alpha_3+ n_4 \alpha_4: n_i \in \Z$ and $n_1+n_3+n_4$ is even$\}$.
\item[$d_2$)] $G$  of type $D_\ell\, (\ell \geq 5): \lambda \in 
\{2n_1 \alpha_1+2n_2
\alpha_2+ \cdots+ 2n_{\ell-2} \alpha_{\ell-2}+n_{\ell-1}
\alpha_{\ell-1}+n_{\ell} \alpha_\ell, n_i \in \Z$ and  $n_{\ell -1}+ n_\ell$
is even$\}$.
\item[e)] $G$ of type $G_2: \lambda \in \Z 6 \alpha_1+\Z 2 \alpha_2$.
\item[f)] $G $ of type $F_4: \lambda \in  \Z 6 \alpha_1 +\Z 6\alpha_2+
\Z 12\alpha_3+
\Z 12 \alpha_4$.
\item[g)]$G$ of type $E_6: \lambda \in 6\Lambda$.
\item[h)]$G$ of type $E_7: \lambda \in 12\Lambda$.
\item[i)] $G$ of type $E_8: \lambda \in 60 Q$.
\end{enumerate}

\noindent
{\bf Acknowledgement:} We acknowledge partial support from the NSF.

\section{Notation}

Let $G$ be a connected, simply-connected, semisimple algebraic group
over the field $\Bbb C$ of
complex numbers. 
We fix a Borel subgroup $B$ of $G$ and a maximal torus
$T$ contained in $B$. Let $B^-$ be the opposite Borel subgroup of $G$ (i.e., the  Borel subgroup $B^-$ such that $B^-\cap B=T$). We denote the unipotent 
radicals of $B,B^-$ by $U,U^-$ respectively. 
We denote the Lie algebra of any algebraic group by the corresponding Gothic character. In particular, the Lie algebras of $G,B,T$ are denoted by 
$\mathfrak{g}, \mathfrak{b}, \mathfrak{t}$ respectively. Let $\triangle=
\triangle(\mathfrak{g})$ be the set of roots of $\mathfrak{g}$ with respect to
$\mathfrak{t}$ and let  $\triangle^+=
\triangle^+(\mathfrak{g})$ be the set of positive roots (i.e., the set of roots of $ \mathfrak{b}$) with $\Pi=\{\alpha_1, \dots, \alpha_\ell\}$ the set of simple roots. Let $W$ be the Weyl group of $G$ with $\{s_1, \dots, s_\ell\} 
\subset W$ the set of simple reflections corresponding to the simple roots 
$\{\alpha_1, \dots, \alpha_\ell\}$ respectively. Let $Q:=\oplus_{i=1}^\ell \, \Z\alpha_i\subset \mathfrak{t}^*$ be the root lattice and let 
$$\Lambda:=\{\lambda\in  \mathfrak{t}^*: \lambda(\alpha_i^\vee)\in \Z \,
{\rm for \ all}\,i\}$$
be the weight lattice, where $\{\alpha_1^\vee, \dots, \alpha_\ell^\vee\}$ are the simple coroots corresponding to the simple roots $\{\alpha_1, \dots, \alpha_\ell\}$ respectively. Let $\Lambda^+\subset \Lambda$ be the set of dominant weights, i.e.,
$$\Lambda^+:=\{\lambda\in  \mathfrak{t}^*: \lambda(\alpha_i^\vee)\in \Z^+ \,
{\rm for \ all}\,i\},$$
where $\Z^+$ is the set of nonnegative integers. For $\lambda \in \Lambda^+$, 
we denote by $V(\lambda)$ the irreducible representation of $\mathfrak{g}$ 
with highest weight $\lambda$. 

 Let $X$ be the group of characters of $T$. Then, $X$ can be identified with 
$\Lambda$ by taking the derivative of characters. For $\lambda \in \Lambda$, we denote the corresponding character of $T$ by $e^\lambda$. Then, 
 $e^\lambda$ uniquely extends to a character of $B$. 

We also consider a (standard) parabolic subgroup $P\supset B$ of $G$. Define
$$\Pi_P=\{{\rm simple\ roots}\ \alpha_i: -\alpha_i \ {\rm is \ a \ root\ of}\
\mathfrak{p}\},$$
$$\Lambda_P=\{\lambda\in \Lambda: \lambda(\alpha_i^\vee)=0 \ {\rm for}\ 
\alpha_i\in \Pi_P\},\,\,
\Lambda_P^+=\Lambda_P\cap \Lambda^+ ,$$ 
and
$$\Lambda_P^o=\{\lambda\in \Lambda_P: \lambda(\alpha_i^\vee)>0 \ {\rm for 
\ each}\ \alpha_i\in \Pi\setminus \Pi_P\}.$$
In particular, for $P=B$, $\Lambda_B=\Lambda$ and  $\Lambda_B^+=\Lambda^+$.
We recall that $e^\lambda$ extends to a character of $P$ iff 
$\lambda\in \Lambda_P$. For $\lambda\in \Lambda_P$, let ${\cal L}_P(\lambda)
=G\times_P \C_{-\lambda}$ be the homogeneous line bundle on the flag variety $G/P$ associated to the principal $P$-bundle $G\rightarrow G/P$ via the one dimensional representation $ \C_{-\lambda}$ of $P$ given by the character $e^{-\lambda}$. Then,  ${\cal L}_P(\lambda)$ is an ample line bundle iff $\lambda\in 
\Lambda_P^o$. For $g\in G$ and $v\in  \C_{-\lambda}$, we denote the equivalence class of $(g,v)$ in $G\times_P \C_{-\lambda}$ by $[g,v]$. 

\section{Proofs }

 For any set of positive roots $\ubeta =\{\beta_1, \dots, 
\beta_n\}$, let $\Z (\ubeta): = \Z \beta_1 + \cdots + \Z \beta_n \subset Q$ and let $\mathfrak{g}
(\ubeta)$ be the semisimple subalgebra of $\mathfrak{g}$ with roots
$$\triangle (\mathfrak{g} (\ubeta)): = \triangle (\mathfrak{g}) \cap \Z
(\ubeta).$$

We recall the  following well-known result.
\begin{lemma} For any $\lambda \in \Lambda^+$ and  any set of positive roots
$\ubeta$  such that the index  of $ \Z (\ubeta)$ in the root 
lattice     $Q$ is finite,  the $0$-weight space of the  submodule $U
(\mathfrak g(\ubeta))\cdot v^+_\lambda \subset V (\lambda)$ is nonzero if and only if
$\lambda \in  \Z (\ubeta)$, where $v^+_\lambda $ is a nonzero highest
weight vector of $V(\lambda)$.

Observe that since $\Z (\ubeta)$ is of finite index in $Q$, the full cartan
subalgebra  $\mathfrak{t} \subset  \mathfrak{g} (\ubeta)$.
\end{lemma}

Let $P$ be a standard parabolic subgroup.
  For the action of $T$ on $G/P$ via the left
multiplication, the isotropy subgroup $I_{gP} \subset T$ at any $gP\in G/P$ 
 is clearly  given by
$$I_{gP} = T \cap  g P g^{-1}.$$

\begin{lemma}  For $gP \in G/P$ and any $\mu \in \Lambda_P$, the isotropy  subgroup  
$I_{gP}$ acts trivially on the fiber ${\cal L}_P (\mu)_{\, \mid_{gP}}$ iff
$${e^{\mu}}_{\, \mid_{P \cap g^{-1} T g}} \equiv 1.$$
\end{lemma}

\noindent
{\it Proof.}  For $t \in I_{gP}= T \cap\, gP g^{-1}$ and any nonzero vector
$v_{-\mu}\in \Bbb C_{-\mu}$,
$$
\begin{array}{lll}
t [g, {v}_{- \mu}] &=& [tg, v_{- \mu}] \\
&=& [g g^{-1} tg, v_{- \mu}] \\
&=& [g, (g^{-1} tg) v_{- \mu}], \,\, {\rm since} \ g^{-1} t g \in P \\
&=& [g, e^{- \mu}  (g^{-1} tg)  v_{- \mu}].
\end{array}
$$
Thus,
$$t [g, v_{- \mu}] = [g, v_{- \mu}]  \,\, {\rm for} \  t \in T \cap g
Pg^{-1} \ {\rm iff} \,\, e^{- \mu} (g^{-1}tg) =1.$$
 This proves the lemma. $\hfill{\Box}$

\begin{lemma}  For any $g=\bar{w} up$, for 
$\bar{w} \in N(T), u \in U^-_{P}$ and $p \in P$, we have
$$P \cap g^{-1} T g = p^{-1} (T \cap u^{-1} Tu) p, \leqno{(1)}
$$
where $N(T)$ is the normalizer of $T$ in $G$ and $U^-_P:=U^-\cap P$.

Moreover, in the above if we just assume that $u\in U^-$, then
$$ P \cap g^{-1} T g \supset p^{-1} (T \cap u^{-1} Tu) p. \leqno{(2)}
$$
\end{lemma}

\noindent
{\it Proof.} We first prove (1). We have
$$
P \cap g^{-1} T g = P \cap p^{-1} u^{-1} \bar{w}^{-1} T \bar{w} up, $$
i.e.,
$$
 P \cap g^{-1} T g= p^{-1} (P \cap u^{-1} Tu)p. \leqno{(3)} 
$$
Now, for any $t \in T,$ if  $u^{-1} tu \in P$, then $u^{-1} tu =t$.  
To prove this, observe that $u^{-1} tut^{-1} \in P \cap U_P^-= (1)$.  Thus,
$p^{-1} (P \cap u^{-1} Tu) p  \subset p^{-1} Tp$ and hence $P \cap g^{-1}Tg \subset p^{-1} Tp \cap  g^{-1} Tg$.  The reverse inclusion is, of course, 
obvious.

The inclusion (2) of course follows from (3).
$\hfill{\Box}$
\vskip2ex

For $Y=G/P$ and an ample line bundle ${\cal L}_P (\lambda)$ on $Y$ (i.e., $\lambda \in \Lambda_P^o$), by 
$Y^{ss} (\lambda)$, we mean the set of semistable points in $Y$ with respect 
to the action of $T$ via the left multiplication on  $Y$ and $T$-linearlized 
ample line bundle ${\cal L}_P (\lambda)$.  Then, as is well known, 
$Y^{ss} (\lambda) \subset Y$ is a Zariski open subset.  
Moreover, by Lemma 3.1, for any $n \geq 1$, $V(n \lambda)^T \neq 0$ iff
 $n \lambda \in Q$. Since $Q$ is of finite index in $\Lambda$, 
we get from the next Lemma 3.4 that
$$Y^{ss} (\lambda) \neq \emptyset.$$

\begin{lemma}  Let $\lambda \in \Lambda_P^o$. Then,  
$gP \in Y^{ss} (\lambda)$ iff 
$gv^+_{n \lambda}$ has a nonzero component in the zero weight space for some 
$n \geq 1$, where $v^+_{n \lambda}$ is a highest weight vector of 
$V(n \lambda)$.
\end{lemma}

\noindent
{\it Proof.}  By definition, the point $ g P  \in Y$ is semistable iff there exists a 
section  $\sigma \in H^0 (Y, {\cal L}_P (n \lambda))^T$ (for some
$n \geq 1)$, such that $\sigma (gP) \neq 0$.  Consider the isomorphism
$$\chi: V (n \lambda)^*  \tilde{\to} H^0 (G/P, {\cal L}_P (n  \lambda)),$$
where $\chi (f) (g P) =[g, (g^{-1}  f)_{ \mid \Bbb C v^+_{n \lambda}}]$.
Thus, 
$$\chi (f) (gP) \neq 0 \Leftrightarrow f (gv^+_{n\lambda}) \neq 0,$$
 and hence
$gP$ is semistable iff $gv^+_{n \lambda}$ has a nonzero component in
the zero weight space. $\hfill{\Box}$

\begin{lemma}  For $gP \in Y^{ss} (\lambda)$,
$$ {e^{n \lambda}}_{ \mid P \cap g^{-1}Tg}  \equiv 1, \ {\rm for \ some} \ n 
\geq 1. \leqno{(1)}$$
In particular,
$${e^{\lambda}}_{\, \mid(P \cap g^{-1} Tg)^o} \equiv 1, 
\leqno{(2)}$$
  where $(P \cap g^{-1} Tg)^o$ is the identity component of $P \cap g^{-1} Tg$.
\end{lemma}

\noindent
{\it Proof.}  By Lemma 3.4, $gP$ is semistable iff for some $n\geq 1$, $[gv^+_{n \lambda}]_0 
\neq 0$, where $[g v^+_{n \lambda}]_0$ denotes its component in the zero
weight space.  For $t \in T$,
$$[tg v^+_{n \lambda} ]_0 = [gv^+_{n \lambda}]_0. \leqno{(3)}$$
But, for $t \in T \cap gP g^{-1}$,
$$
[tg v^+_{n \lambda}]_0 = [g g^{-1} tg v^+_{n \lambda}]_0 
= e^{n \lambda} (g^{-1} tg) [gv^+_{n \lambda}]_0. \leqno{(4)}
$$
Combining (3) and (4), we get $e^{n \lambda}  (g^{-1} t g)=1$.  This
proves (1).  

The identity (2) follows immediately from (1) since a connected 
torus is a divisible group. 
$~~~~\hfill{\Box}$

\begin{lemma}  For any $x \in \mathfrak{u}^-$, let  $\ubeta_x$ be the subset of $\triangle^+ $ such that
$x = \sum_{\alpha \in \ubeta_x} x_{-\alpha}$, where $x_{- \alpha} \in \mathfrak{g}_{- \alpha}$ is nonzero.
Then, for $u = $ Exp $x$,
$$T \cap u^{-1} T u = \displaystyle{\cap_{ \alpha \in \ubeta_x}} \{ t  \in T: e^\alpha (t) =1\}. \leqno{(1)}
$$
In particular, for a weight $\mu \in \Lambda, {e^{\mu}}_{\, \mid  T \cap u^{-1} 
Tu}
\equiv 1$ iff $\mu \in \Z \ubeta_x$.
\end{lemma}

\noindent
{\it Proof.}  Take $t \in T$ such that $utu^{-1} \in T$.  Then,
$$utu^{-1}t^{-1} \in T \cap U^- = \{1\}.$$  Thus,
$utu^{-1} =t$,  i.e., $tut^{-1} =u.$
But,  
 \begin{align*}tu t^{-1} &= {\rm Exp} \ ({\rm Ad}\, t \cdot x) \\
&= {\rm Exp} \ ({\sum_{\alpha \in \ubeta_x}}e^{- \alpha} (t) 
x_{- \alpha}) \\
&={\rm Exp} (x). 
\end{align*}
Since Exp$_{\, \mid \mathfrak{u}^-} $ is a bijection,
$${\sum_{\alpha \in \ubeta_x}} e^{ - \alpha} (t)  x_{- \alpha} = 
{\sum_{\alpha \in \ubeta_x}}x_{-\alpha}.$$
Thus, $e^{-\alpha} (t) =1$ for any $\alpha \in \ubeta_x$.  This proves the inclusion
$$T \cap u^{-1} Tu \subset \displaystyle{\cap_{\alpha \in \ubeta_x} }\{t \in T:
e^{\alpha} (t) =1\}.$$
The reverse inclusion follows by reversing the above calculation.  

To prove 
the `In particular' statement, consider the isomorphism $\xi:T\simeq 
{\rm Hom}_{\Bbb Z}(\Lambda, \Bbb C^*)$,  $\xi (t) \mu =e^\mu (t)$, for
  $t \in T $ and $\mu \in \Lambda$. From this identification and (1) it is
easy to see that $\xi (T \cap u^{-1} Tu)  = 
\Hom_{\Z} (\frac{\Lambda}{\Z \ubeta_x},
\C^*)$.  Now, if $\mu \in  \Lambda \setminus  (\Z \ubeta_x)$, there exists a homomorphism $f: \frac{\Lambda}{\Z \ubeta_x} \to \C^*$  Such that $f (\mu) \neq 
1$. From this we conclude the `In particular' statement. $\hfill{\Box}$

\vskip2ex

Since $Y^{ss} (\lambda)$ is Zariski open in $Y$, we can find a semistable point of the form Exp$(x)P$, where $x \in \mathfrak{u}^-$ and $\ubeta_x = \triangle^+$.

More generally, take any subset $\ubeta \subset \triangle^+$
such that $\Z (\ubeta)$ is of finite index in $Q$.  Then, for some $n\geq 1, n \lambda
\in \Z (\ubeta)$ and hence by Lemma 3.1,
$$\bigl(U (\mathfrak{g} (\ubeta))\cdot v^+_{n \lambda}\bigr)^T \neq (0).$$
In particular,
$$\bigl(G(\ubeta) / B(\ubeta)\bigr) \cap Y^{ss} (\lambda) \neq \emptyset,$$
where $G(\ubeta)$ is the connected (semisimple) subgroup of $G$ with Lie 
algebra $\mathfrak{g} (\ubeta) $ and $B (\ubeta):= B \cap G(\ubeta)$
is a Borel subgroup of $G(\ubeta)$.  Again, by Zariski density, we can find an element of the form Exp$(x)P \in Y^{ss} (\lambda)$ such that $x\in 
\mathfrak{u}^-$ and $\ubeta_x$ is the set of all the positive roots of $G(\ubeta)$ (i.e., all the roots of $B (\ubeta)$).

\begin{lemma}  For any subset $S \subset \triangle^+$, the quotient group
$$T_S/ T_S^o \simeq \Tor (\Lambda / \Z S),$$
where $ T_S: = \displaystyle{\cap_{\alpha \in S}}  \{ t
\in T: e^\alpha (t) =1\}$, $T_S^o$ denotes its identity component and Tor denotes the torsion subgroup.
\end{lemma}

\noindent
{\it Proof.} Recall the identification $\xi: T \tilde{\to} \Hom_\Z
(\Lambda, \C^*)$ from the proof of Lemma 3.6.   Under this identification
$$\xi (T_S) = \Hom_{\Z} \bigl(\frac{\Lambda}{\Z S} , \C^*\bigr).$$
Decompose
$$\frac{\Lambda}{\Z S} \simeq \Tor \bigl(\frac{\Lambda}{\Z S} \bigr)\oplus F,$$
where $F$ is a free $\Z$-module.  Thus,
$$\Hom_{\Z} (\frac{\Lambda}{\Z S}, \C^*) \simeq \Hom_{\Z} (\Tor (\frac{\Lambda}{\Z S}), \C^*) \times \Hom_{\Z} (F, \C^*).$$
  
But, $\Hom_{\Z} (F, \C^*) \simeq (\C^*)^m$, where $m:=  {\rm rk}\,F$ and, of course, 
$$\Hom_{\Z} (\Tor (\frac{\Lambda}{\Z S}), \C^*) \simeq \Tor (\frac{\Lambda}{\Z S}).$$
Thus, $T_S^o \simeq \Hom_\Z (F,\C^*)$ and hence
$$T_S / T_S^o \simeq \Tor (\Lambda / \Z S). \eqno{\Box}$$

Let $H$ be a reductive group, $X$ a projective variety and ${\cal L}$ an 
ample $H$-equivariant line bundle on $X$. Then, recall that the {\it GIT quotient} $X//H$ is by definition the uniform categorical quotient of the (open)
set of semistable points $X^{ss}$ by $H$ (cf. [MFK, Theorem 1.10]).

We recall the following `descent' lemma of Kempf (cf. [DN, Theorem 2.3]).

\begin{lemma} Let $X, H $ and $ {\cal L}$ be as above. Then, an $H$-equivariant vector bundle ${\cal S}$ on $X$ descends to a vector bundle on $X//H$ 
(i.e., there exists a vector bundle ${\cal S'}$ on $X//H$ such that its pull-back to $X^{ss}$ under the canonical $H$-equivariant structure is 
$H$-equivariantly isomorphic to the restriction of ${\cal S}$ to $X^{ss}$) 
iff for any $x\in X^{ss}$, the isotropy subgroup $I_x$ acts trivially on the fiber ${\cal S}_x$. 

In fact (though we do not need), for the `if' part, it suffices to assume that $I_x$ acts trivially for only those $x\in X^{ss}$ such that the orbit $H\cdot x$ is closed in $X^{ss}$.
\end{lemma}
Let $P$ be a standard parabolic subgroup and let 
${\cal L}_P (\lambda)$ be a homogeneous ample line bundle (i.e., $\lambda \in \Lambda^o_P$) on $Y=G/P$.  Denote the GIT quotient of $Y$ by $T$ with respect 
to ${\cal L}_P (\lambda)$ by $Y(\lambda)//T$. The following is one of our main results.

 \begin{theorem}  With the notation as above, the line bundle 
${\cal L}_P (\lambda)$
descends to a line bundle on  $Y (\lambda) // T $  iff
 for all the semisimple subalgebras $\mathfrak{s}$ of $\mathfrak{g}$ containing $\mathfrak{t}$ (in particular,  rank $\mathfrak{s}$ = rank
 $\mathfrak{g})$,
$$\lambda \in \Z\triangle^+ (\mathfrak{s}),$$
where $\triangle^+ (\mathfrak{s}):=\triangle^+ \cap \triangle(\mathfrak{s})$ 
is the set of positive roots of $\mathfrak{s}$.
\end{theorem}

\paragraph*{Proof.}  By Lemma 3.8, ${\cal L}_P (\lambda)$
descends to $Y (\lambda) // T$ iff for all the semistable points $gP=
\bar{w} u P \in G/P$ (for $\bar{w} \in N(T)$ and $u \in U^-_P)$, 
the isotropy subgroup $I_{gP}$ acts trivially on the fiber ${\cal L}_P 
(\lambda)_{\mid gP}$. By Lemmas 3.2 and 3.3, this is equivalent to the 
requirement that
${e^\lambda}_{\, \mid T \cap u^{-1} Tu} \equiv 1$.  Further, by Lemma 3.6, this
is equivalent to the requirement that
$ \lambda \in \Z \ubeta_x,$
where $x \in \mathfrak{u}^-$ is the element with Exp $x=u$.

By the discussion above Lemma 3.7, for any semisimple subalgebra $\mathfrak{s}$ 
of $\mathfrak{g}$ containing $\mathfrak{t}$, there exists an element $x \in \mathfrak{u}^-$
such that Exp$(x)P \in G/P$ is semistable and, moreover, $\ubeta_x = \triangle^+
(\mathfrak{s})$.  Thus, we get by Lemmas 3.2, 3.3 and 3.6 that $\lambda \in 
\Z\triangle^+ (\mathfrak{s})$ for any such $\mathfrak{s}$
if the line bundle ${\cal L}_P (\lambda)$ descends to $Y (\lambda) //T$.

Conversely, take any semistable point  $\bar{w}$Exp$(x)P \in G/P$ for 
$\bar{w}\in N(T)$ and $x\in \mathfrak{u}^-_P$. If $\Z \ubeta_x
 \subset Q$ is not of finite index in $Q$, choose simple roots $\ualpha_x= \{\alpha_{i_1}, \dots, \alpha_{i_j}\}$ such that $\Q \ubeta_x
\cap \Q \ualpha_x = (0)$ and $\Q \ualpha_x+\Q \ubeta_x= {\oplus_{i=1}^{\ell}}\,\Q \alpha_i$, where $\Q \ubeta_x:= {\sum_{\beta \in \ubeta_x}} \Q \beta 
\subset \mathfrak{t}^*$ and $ \Q\ualpha_x:= {\oplus_{n=1}^{j}} \Q 
\alpha_{i_n} \subset \mathfrak{t}^*$. With this choice of $\ualpha_x$, the 
torsion submodule
$$\Tor (\Lambda / \Z \ubeta_x) \hookrightarrow \Tor (\Lambda / (\Z \ubeta_x+\Z \ualpha_x)). 
\leqno{(1)}$$
 To see this observe that we have the following short exact sequence:
$$0\rightarrow \Z \ualpha_x \rightarrow \Lambda / \Z \ubeta_x \to \Lambda / 
(\Z \ualpha_x + \Z \ubeta_x)\rightarrow 0.$$
From this the assertion (1) follows easily.

Let $\mathfrak{s}$ be the semisimple  subalgebra of   $\mathfrak{g}$ containing  
$\mathfrak{t}$  such that
$$\Z\triangle^+ (\mathfrak{s})= 
 \Z \ubeta_x + \Z \ualpha_x.$$
Choose an element $y= y_{\mathfrak{s}} \in \mathfrak{u}^-$ such that
$\ubeta_y = \triangle^+ (\mathfrak{s})$ and (Exp $y)P \in Y^{ss} (\lambda)$. With 
this choice of $y$, by Lemma 3.6, we see that $T \cap v^{-1} Tv \subset T \cap
u^{-1} T u$, where $u: = $ Exp $x$ and $v:= $ Exp $y$.  We
further show that the inclusion $T \cap v^{-1} Tv \subset T \cap u^{-1}
Tu$ induces a  splitting:
$$(T \cap u^{-1}T u) / (T \cap u^{-1} Tu)^o \hookrightarrow (T \cap v^{-1}
Tv) / (T\cap v^{-1} Tv)^o. \leqno{(2)}$$
By Lemma 3.7,
$$(T \cap v^{-1} Tv) /  (T \cap v^{-1} Tv)^o   \simeq \Tor (\frac{\Lambda}{\Z \ubeta_y})$$
and 
$$ (T\cap u^{-1} Tu) / (T \cap u^{-1} Tu)^o \cong \Tor (\frac{\Lambda}{\Z \ubeta_x}).$$
Combining the above identifications with (1), we get (2).

Since $\bar{w}u P$ is a semistable point, by Lemmas 3.3 and 3.5, 
${e^\lambda}_{\vert (T \cap
u^{-1} Tu)^o} =1$.  Moreover, by the assumption and Lemma 3.6, 
${e^{\lambda}}_{\, \mid T \cap 
v^{-1} Tv} \equiv 1$.  Thus, ${e^\lambda}_{\, \mid T \cap u^{-1} Tu} \equiv 1$
by (2).
This proves the theorem. $\hfill{\Box}$

\vskip1ex

Using the above theorem, we explicitly get exactly which
 ${\cal L}_P(\lambda)$  descends to the GIT quotient $Y (\lambda) // T$.

In the following, we follow the indexing
convention as in Bourbaki [B, Planche I - IX].

\begin{theorem}  Let $G$ be a connected, simply-connected simple algebraic group, 
$P \subset G$ a standard parabolic subgroup and let ${\cal L}_P (\lambda)$
be a homogeneous ample line bundle on the flag variety $Y= G/P$ (i.e., $\lambda \in 
\Lambda_P^o$).  Then, the
line bundle ${\cal L}_P (\lambda)$ descends to a line bundle on the GIT quotient
$Y(\lambda) //T$ if and only if $\lambda$ is of the following form depending upon the 
type of $G$ (in addition to $\lambda \in \Lambda_P^o$).

\begin{enumerate}
\item[a)] $G$ of type $A_{\ell}\, (\ell \geq 1): \lambda \in Q$.
\item[b)]$G $ of type $B_\ell \,(\ell \geq 3): \lambda \in 2Q$.
\item[c)]$G $ of type $C_\ell \,(\ell \geq 2): \lambda \in \Z 2\alpha_1+
\cdots +\Z 2 \alpha_{\ell-1}+ \Z \alpha_\ell = 2\Lambda $.
\item[$d_1$)]  $G $ of type $D_{4}: \lambda \in \{n_1 \alpha_1+ 2n_2 \alpha_2+
n_3 \alpha_3+ n_4 \alpha_4: n_i \in \Z$ and $n_1+n_3+n_4$ is even$\}$.
\item[$d_2$)] $G$  of type $D_\ell\, (\ell \geq 5): \lambda \in 
\{2n_1 \alpha_1+2n_2
\alpha_2+ \cdots+ 2n_{\ell-2} \alpha_{\ell-2}+n_{\ell-1}
\alpha_{\ell-1}+n_{\ell} \alpha_\ell, n_i \in \Z$ and  $n_{\ell -1}+ n_\ell$
is even$\}$.
\item[e)] $G$ of type $G_2: \lambda \in \Z 6 \alpha_1+\Z 2 \alpha_2$.
\item[f)] $G $ of type $F_4: \lambda \in  \Z 6 \alpha_1 +\Z 6\alpha_2+
\Z 12\alpha_3+
\Z 12 \alpha_4$.
\item[g)]$G$ of type $E_6: \lambda \in 6\Lambda$.
\item[h)]$G$ of type $E_7: \lambda \in 12\Lambda$.
\item[i)] $G$ of type $E_8: \lambda \in 60 Q$.
\end{enumerate}
\end{theorem}

\paragraph*{Proof.} The theorem follows by using Theorem 3.9  and the classification of
the semisimple subalgebras of  $\mathfrak{g}$ of maximal rank due to  Borel-Siebenthal 
(cf. [W, Theorem
8.10.8]).  In the following chart, we only list all the proper maximal semisimple
subalgebras $\mathfrak{s}$ of $\mathfrak{g}$ containing
the (fixed) Cartan subalgebra $\mathfrak{t}$ up to a conjugation under the Weyl
group.  In the following, $\theta$ denotes the highest
root and $\mathfrak{s}_i$ denotes the  semisimple subalgebra of
$\mathfrak{g}$ containing $\mathfrak{t}$ with simple roots $\{\alpha_1, \dots,
\widehat{\alpha}_i, \dots, \alpha_\ell, -\theta \}$.

\vskip1ex

a) $A_\ell \,\,(\ell \geq 1)$: none 

b)$B_\ell\,\, (\ell \geq 3): \mathfrak{s}_i, 2 \leq i \leq \ell$

c) $C_\ell\,\, (\ell \geq 2): \mathfrak{s}_i, 1\leq i \leq \ell-1$

d)$D_\ell \,\, (\ell \geq 4): \mathfrak{s}_i, 2 \leq i \leq \ell -2$

e)$G_2 : \mathfrak{s}_i$, i=1,2

f) $F_4: \mathfrak{s}_i,  i=1,2,4$

g) $E_6: \mathfrak{s}_i,  i=2,3,4,5$

h) $E_7 : \mathfrak{s}_i,  i=1,2,3,5,6$

i) $E_8: \mathfrak{s}_i, i=1,2,5,7,8$.

\vskip1ex

We denote by $L(\mathfrak{g}; \alpha_1, \dots, \alpha_\ell)$ the intersection of
$\Z\triangle^+ (\mathfrak{s})$ 
(inside the root lattice $Q)$ as $\mathfrak{s}$ varies over all
possible semisimple subalgebras $\mathfrak{s}$ of $\mathfrak{g}$ 
 containing the (fixed) cartan subalgebra $\mathfrak{t}$. 

In the following, we  determine $L(\mathfrak{g}; 
\alpha_1, \dots, \alpha_\ell)$ for each simple $\mathfrak{g}$. 

\vskip1ex

\noindent
{\bf a) $A_\ell \,\,(\ell \geq 1):$} In this case 
$L (A_\ell; \alpha_1, \dots, \alpha_\ell)= Q$.

\vskip1ex

\noindent
{\bf c) $C_\ell \,\,(\ell \geq 2):$} By the above chart,
$$
\begin{array}{lll}
L(C_{\ell}; \alpha_1, \dots, \alpha_\ell) &=& [(\Z \theta+L (C_{\ell -1};
\alpha_2, \dots, \alpha_\ell)) \cap \\
&&(L(C_2; \alpha_1, - \theta)+ L (C_{\ell-2}; \alpha_3, \dots, \alpha_\ell)) \cap \cdots \\
&&(L (C_{\ell-2}; \alpha_{\ell-3}, \dots, \alpha_1, - \theta)+L (C_2; \alpha_{\ell-1}, \alpha_\ell)) \cap \\
&&(L(C_{\ell-1};  \alpha_{\ell -2}, \dots, \alpha_1, - \theta)+ \Z \alpha_\ell)]_W,
\end{array}
$$
where $[M]_W$ denotes $\displaystyle{\cap_{w \in W}} \,w M$. By induction, for $j< \ell$,
$$
L(C_j; \alpha_{\ell-j+1}, \dots, \alpha_\ell) = \Z 2 \alpha_{\ell -j+1}+ \cdots+ \Z 2 \alpha_{\ell-1}+ \Z\alpha_\ell $$
and 
\begin{align*} L (C_j; \alpha_{j-1}, \dots, \alpha_1, - \theta) &= \Z 2 
\alpha_{j-1}+ \cdots + \Z 2 \alpha_1+\Z \theta \\
&= \Z 2 \alpha_{j-1}+ \cdots+ \Z 2 \alpha_1+\Z (2 \alpha_j+ \cdots + 2 
\alpha_{\ell-1}+ \alpha_\ell), 
\end{align*}
since 
$$\theta (C_\ell; \alpha_1, \dots, \alpha_\ell) = 2 \alpha_1+
\cdots + 2 \alpha_{\ell-1} + \alpha_\ell.
$$
Thus,
$$L (C_\ell; \alpha_1, \dots, \alpha_\ell) = [\Z 2 \alpha_1+ \cdots + \Z 2 \alpha_{\ell -1} + \Z \alpha_\ell]_W.
$$
But, $\Z 2 \alpha_1+ \cdots + \Z 2 \alpha_{\ell -1} +\Z \alpha_\ell$ is $W$-stable and hence 
$$L (C_\ell; \alpha_1, \dots, \alpha_\ell )= \Z  2 \alpha_1+
\cdots + \Z 2 \alpha_{\ell -1}+\Z\alpha_\ell.$$
This proves the (c)-part of the theorem.

\vskip1ex

\noindent
{\bf d) ${D_\ell \,\,(\ell \geq 4)}:$} We first consider the case of $D_4$. In 
this case
$$
\begin{array}{lll}
L(D_4; \alpha_1, \dots, \alpha_4)&=& [L (A_1; \alpha_1)+ L(A_1;-\theta)+L(A_1; \alpha_3)+
L(A_1; \alpha_4)]_W \\
&=& [\Z \alpha_1+ \Z \theta +\Z \alpha_3+ \Z \alpha_4]_W \\
&=& [\Z \alpha_1+ \Z 2 \alpha_2+ \Z \alpha_3+\Z \alpha_4]_W.
\end{array}
$$
Now, the  sublattice $L':= \{n_1 \alpha_1 +n_2 \alpha_2+ n_3 \alpha_3+ n_4 \alpha_4: n_i 
\in \Z, n_2$ is even and ${\sum_{i=1}^{4}}n_i$ is even$\}
\subset \Z \alpha_1+ \Z 2 \alpha_2+\Z \alpha_3+ \Z \alpha_4$ is $W$-stable
(as is easy to see).  Moreover, clearly the index of $L'$ in $\Z \alpha_1+
\Z 2 \alpha_{2}+ \Z \alpha_3+\Z \alpha_4$ is 2 and $\Z \alpha_1+\Z 2 \alpha_2+
\Z \alpha_3+ \Z \alpha_4$ is not $W$-stable (as can be seen by applying the second simple 
reflection $s_2)$.  Thus $[\Z \alpha_1+\Z 2 \alpha_2 + \Z \alpha_3+
\Z \alpha_4]_W = L'$, proving the theorem in this case.

We now come to the case of general $D_\ell$.  By  the above chart and the 
$A_\ell$ case,
\begin{align*} L(D_\ell; \alpha_1, \dots, \alpha_\ell)&= [(L (D_2; \alpha_1, -\theta)+ 
L(D_{\ell-2}; 
\alpha_3, \dots, \alpha_\ell)) \cap \cdots\\
&\cap (L(D_{\ell-2}; \alpha_{\ell-3}, \dots, \alpha_1, - \theta)+ L(D_2; \alpha_{\ell-1}, 
\alpha_\ell))]_W,
\end{align*}
where $D_k$ for $k=2,3$ is interpreted as $A_k$. Set
$$ L' (\alpha_{\ell -3}, \alpha_{\ell -2}, \alpha_{\ell -1},\alpha_\ell)=
\{\sum_{i=\ell-3}^\ell\,n_i\alpha_i: n_i\in \Z, \sum n_i \,\,
{\rm is\, even \,and}\, n_{\ell -2} \,\,{\rm is\, even}\},$$
and if $\ell -k >4$, set
$$ L'(\alpha_{k+1}, \dots,\alpha_\ell)=
\{\sum_{i=k+1}^\ell\,n_i\alpha_i: n_i\in \Z, \sum n_i \,\,
{\rm is \, even \, and}\, n_{k+1}, \dots, n_{\ell -2} \,\,{\rm are \,even}
\},$$
and if $\ell -k <4$, set
$$ L'(\alpha_{k+1}, \dots,\alpha_\ell)= \Z\alpha_{k+1}+ \cdots +\Z \alpha_\ell.$$

 By induction and the $A_{\ell-k}$ case, for $2\leq 
\ell - k < \ell$, 
$$L(D_{\ell-k}; \alpha_{k+1}, \dots,\alpha_\ell) = 
L'(\alpha_{k+1}, \dots,\alpha_\ell).$$
Thus,
$$
\begin{array}{lll}
L(D_\ell; \alpha_1, \dots, \alpha_\ell) &=& \bigl[(\Z \alpha_1+\Z \theta + L' 
(\alpha_3, 
\dots, \alpha_\ell)) \cap \\
&&(\Z \alpha_1+ \Z \alpha_2+ \Z \theta+ L' (\alpha_4, \dots, \alpha_\ell)) \cap \\
&&(L'(\alpha_3, \alpha_2, \alpha_1, -\theta)+ L' ( \alpha_5, \dots , 
\alpha_\ell)) \cap \cdots\\
&&\cap (L' (\alpha_{\ell -5}, \dots, \alpha_1, - \theta) + 
L' (\alpha_{\ell-3}, \alpha_{\ell - 2}, \alpha_{\ell -1},\alpha_\ell))
\\
&& \cap (L' (\alpha_{\ell -4},  \dots, \alpha_1, - \theta)+  
\Z\alpha_{\ell-2}+ \Z \alpha_{\ell-1}+\Z \alpha_\ell) \\
&& \cap (L'(\alpha_{\ell-3}, \dots,  \alpha_1,- \theta)+ 
\Z\alpha_{\ell-1}+ \Z\alpha_\ell)]_W \\
&=& [L'(\alpha_1, \dots, \alpha_\ell)]_W.
\end{array}
$$
It is easy to see that $L'(\alpha_1, \dots, \alpha_\ell)$ is $W$-invariant and hence
$$[L'(\alpha_1, \dots, \alpha_\ell)]_W =L' (\alpha_1, \dots, \alpha_\ell).$$
This proves the theorem in the case of $D_\ell$.

\vskip1ex

\noindent
{\bf (b) ${B_\ell\,\, (\ell \geq 3)}$:}  By the above chart, we get that
$$
\begin{array}{lll}
L(B_\ell; \alpha_1, \dots, \alpha_\ell) &=& [(\Z \alpha_1+ \Z \theta+ L 
(B_{\ell-2}; \alpha_3, \dots, \alpha_\ell)) \cap \\
&&(L (A_3; \alpha_1, \alpha_2, - \theta)+L (B_{\ell-3}; \alpha_4, \dots,
\alpha_\ell)) \cap \\
&&(L(D_4; \alpha_3, \alpha_2, \alpha_1, - \theta) + L (B_{\ell -4}; \alpha_5, \dots, 
\alpha_\ell)) \cap \cdots \\
&&(L(D_{\ell-3}; \alpha_{\ell -4}, \dots, \alpha_1, - \theta)+L (B_3; \alpha_{\ell -2}, \alpha_{\ell -1}, \alpha_{\ell} )) \cap \\
&&(L(D_{\ell -2}; \alpha_{\ell -3}, \dots, \alpha_1, - \theta)+ L (C_2; \alpha_\ell, \alpha_{\ell-1} )) \cap \\
&&(L(D_{\ell-1}; \alpha_{\ell -2}, \dots, \alpha_1, - \theta)+ \Z \alpha_\ell) \cap \\
&&L (D_\ell; \alpha_{\ell-1}, \dots, \alpha_1, - \theta)]_W. \\
\end{array}
$$

By using the result for $D_\ell$ and $A_3$  and also, by  induction, for
$B_j (j < \ell)$, we get that (with $L'$ as in the proof of the $D_\ell$ case)
$$
\begin{array}{lll}
L(B_\ell ; \alpha_1, \dots, \alpha_\ell)& =&[ (\Z \alpha_1+ \Z \theta+ \Z 2 \alpha_3+ \cdots + \Z 2 \alpha_\ell) \cap \\
&&(\Z \alpha_1+ \Z \alpha_2+ \Z \theta + \Z 2 \alpha_4+ \cdots +\Z 2 \alpha_\ell) \cap \\
&&(L' (\alpha_3, \alpha_2, \alpha_1, - \theta)+ \Z 2 \alpha_5+ \cdots +
\Z 2 \alpha_\ell) \cap \cdots \\
&&(L' (\alpha_{\ell -4}, \dots , \alpha_1, - \theta) + \Z 2 \alpha_{\ell -2}+
\Z 2 \alpha_{\ell-1} + \Z 2 \alpha_\ell) \cap \\
&&(L' (\alpha_{\ell-3}, \dots, \alpha_1 - \theta) + \Z 2 \alpha_\ell +\Z
\alpha_{\ell -1})\cap \\
&&(L'(\alpha_{\ell-2}, \dots, \alpha_1, - \theta) +\Z \alpha_\ell ) \cap 
 L' (\alpha_{\ell -1 }, \dots,  \alpha_1, - \theta)]_W \\
&=& [\Z 2 \alpha_1+\cdots + \Z 2 \alpha_\ell]_W \\
&=& 2 [Q]_W \\
&=& 2Q ,\ {\rm since} \,Q \,{\rm is} \ W{\rm -stable}. 
\end{array}
$$ 
This proves the theorem for the case of $B_\ell$.

\paragraph*{(e)  $G_2:$}  By the chart and the theorem for the case of $A_2$,
$$
\begin{array}{lll}
L(G_2; \alpha_1, \alpha_2) &=& [(L (A_2; \alpha_2, - \theta)) \cap (\Z \alpha_1+ \Z 
\theta)]_W \\
&=& [(\Z \alpha_2+\Z \theta) \cap (\Z \alpha_1+ \Z \theta)]_W \\
&=& [(\Z \alpha_2 + \Z 3 \alpha_1) \cap (\Z  \alpha_1 + \Z 2 \alpha_2)]_W \\
&=& [\Z 3 \alpha_1+ \Z 2 \alpha_2 ]_W \\
&=& \Z 6 \alpha_1 + \Z 2 \alpha_2,
\end{array}
$$
since $\Z 6 \alpha_1+ \Z 2 \alpha_2$ is $W$-stable, whereas $\Z 3 \alpha_1+ \Z
2 \alpha_2$ is not $W$-stable and $\Z 6 \alpha_1+ \Z 2 \alpha_2$ is of index 2 inside $\Z 3 \alpha_1+ \Z 2 \alpha_2$.

\paragraph*{(f) $ F_4:$} By the chart,
$$
\begin{array}{lll}
L(F_4; \alpha_1, \alpha_2, \alpha_3, \alpha_4) &=& [(\Z \theta +L( C_3; \alpha_4, 
\alpha_3, \alpha_2)) \cap \\
&&(L (A_2; - \theta, \alpha_1)+ L (A_2; \alpha_3, \alpha_4)) \cap \\
&&L(B_4; - \theta, \alpha_1, \alpha_2, \alpha_3)]_W 
\end{array}
$$
By the theorem for $A_2, C_3$ and $B_4$, we get
$$
\begin{array}{lll}
L (F_4; \alpha_1, \alpha_2, \alpha_3, \alpha_4) &=& [(\Z \theta+ \Z 2 \alpha_4+
\Z 2 \alpha_3+ \Z \alpha_2) \cap \\
&&(\Z \theta+ \Z \alpha_1+ \Z \alpha_3+ \Z \alpha_4) \cap \\
&&(\Z 2 \theta+ \Z 2 \alpha_1+ \Z 2 \alpha_2+ \Z 2 \alpha_3)]_W \\
&=& [(\Z 2 \alpha_1+\Z 2 \alpha_4+ \Z 2 \alpha_3+ \Z \alpha_2) \cap \\
&& (\Z 3 \alpha_2+ \Z \alpha_1+ \Z \alpha_3+ \Z \alpha_4) \cap \\
&&(\Z 4 \alpha_4+\Z 2 \alpha_1+ \Z 2 \alpha_2+ \Z 2 \alpha_3)]_W \\
&=& [M]_W,\ {\rm where} \, M: = \Z 2 \alpha_1+\Z 6 \alpha_2+\Z 2 \alpha_3+\Z 4 \alpha_4.
\end{array}
$$
Now,
$$\Z 6 \alpha_1+ \Z 6 \alpha_2+ \Z 12 \alpha_3+\Z 12 \alpha_4  \subset
[M]_W.$$
Moreover, $\Z 6 \alpha_1+\Z 6 \alpha_2+ \Z 12 \alpha_3+ \Z 12 \alpha_4$ is $W$-stable 
(as can be easily seen).  Conversely, take $\mu \in [M]_W$.  Then,
$w \mu \in M$ for any $w\in W$.  Since 
$$s_i w \mu = w \mu - (w \mu) (\alpha_i^\vee) \alpha_i;$$
for any $\mu \in [M]_W$ and any $w \in W$,
$$w \mu (\alpha_2^{\vee}) \in 6 \Z \ {\rm and} \ w \mu (\alpha^{\vee}_4) \in 
4 \Z.$$
Since $\alpha_1^{\vee} \in W. \alpha_2^{\vee}$ (cf [H, $\S$10.4, Lemma c]),
we get
$$\mu (\alpha_1^\vee) \in  6 \Z,\, \mu (\alpha^\vee_2 ) \in 6  \Z
\,\, {\rm and}  \ \mu (\alpha_4^\vee) \in 4 \Z.$$ 
Take $\mu = n_1 2 \alpha_1+n_2 6 \alpha_2+ n_3 2 \alpha_3+ n_4 4 \alpha_4
\in [M]_W$.  Then, from the above, we get
$$4n_1 -6n_2 \in 6 \Z,\, \,-2n_1+12n_2-2n_3 \in 6 \Z\,\, {\rm and}\,
-2 n_3+8n_4 \in 4 \Z.$$  This gives that $n_1 \in 3 \Z$ and $n_3\in 6 \Z$.  
Thus, $\mu \in \Z 6 \alpha_1+ \Z 6 \alpha_2+ \Z 12 \alpha_3+ \Z 4 \alpha_4$.
Considering $s_3 (\Z 6 \alpha_1+ \Z 6 \alpha_2+ \Z 12 \alpha_3+ \Z 4 \alpha_4)$, we see 
that 
$$[\Z 6 \alpha_1+ \Z 6 \alpha_2+ \Z 12 \alpha_3+ \Z 4 \alpha_4]_W
\subset \Z 6 \alpha_1+ \Z 6 \alpha_2+ \Z 12 \alpha_3+ \Z 12 \alpha_4.$$
Thus, we conclude that 
$$[M]_W=\Z 6 \alpha_1+\Z 6 \alpha_2+ \Z 12 \alpha_3+\Z 12 \alpha_4.$$
This completes the proof of the theorem in the case of $F_4$.

\paragraph*{(g)  $E_6:$} By the chart and the case of $A_\ell$,
$$
\begin{array}{lll}
L(E_6; \alpha_1, \dots, \alpha_6) &=& [(\Z \theta+ L (A_5; \alpha_1, \alpha_3,
\alpha_4, \alpha_5, \alpha_6)) \cap \\
&&(\Z \alpha_1+ L (A_5; - \theta, \alpha_2, \alpha_4, \alpha_5, \alpha_6)) \cap\\
&& (L(A_2; \alpha_1, \alpha_3)+L (A_2; - \theta,  \alpha_2)+ L (A_2; \alpha_5, \alpha_6)) 
\cap \\
&&(L(A_5; \alpha_1, \alpha_3, \alpha_4, \alpha_2, -\theta) + \Z \alpha_6)]_W \\
&=& [(\Z \theta+ \Z \alpha_1+ \Z \alpha_3+ \Z \alpha_4+ \Z \alpha_5+ \Z \alpha_6) \cap \\
&&(\Z \alpha_1+ \Z \theta+ \Z \alpha_2+ \Z \alpha_4+ \Z \alpha_5+ \Z \alpha_6) \cap \\
&& (\Z \alpha_1+ \Z \alpha_3+\Z \theta + \Z \alpha_2+\Z \alpha_5+\Z \alpha_6) \cap \\
&&(\Z \alpha_1+ \Z \alpha_3+ \Z \alpha_4+ \Z \alpha_2+ \Z \theta +\Z \alpha_6) ]_W \\
&=& [(\Z \alpha_1+ \Z 2 \alpha_2+ \Z \alpha_3+\Z \alpha_4+ \Z \alpha_5+ \Z \alpha_6) \cap \\
&&(\Z \alpha_1+ \Z \alpha_2+ \Z 2 \alpha_3 + \Z \alpha_4+ \Z \alpha_5+ \Z
\alpha_6) \cap \\
&&(\Z \alpha_1+ \Z \alpha_2+ \Z \alpha_3+\Z 3 \alpha_4+ \Z \alpha_5+ \Z \alpha_6) \cap \\
&&(\Z \alpha_1+\Z \alpha_2+ \Z \alpha_3+ \Z \alpha_4+ \Z 2 \alpha_5+\Z \alpha_6)]_W \\
&=& [M_6]_W, 
\end{array}
$$
where
$M_6:= \Z \alpha_1+\Z 2 \alpha_2+ \Z 2 \alpha_3+ \Z 3 \alpha_4+ \Z 2 \alpha_5+
\Z \alpha_6.$

Clearly,
$6 \Lambda \subset M_6$  
and since
$\Lambda$ is $W$-stable,
$$6\Lambda \subset [M_6]_W. \leqno{(1)}$$ 
Conversely, take $\mu \in [M_6]_W$. Then, for any $w \in W, w \mu \in M_6$.
Since $s_i (w \mu) = w \mu - (w \mu) (\alpha_i^\vee) \alpha_i$; for any
$\mu \in [M]_W$, and any $w \in W$,
$w \mu (\alpha_2^\vee) \in 2 \Z$ and $w \mu (\alpha^\vee_4) \in 3 \Z$.  
Since $E_6$ is simplylaced, the Weyl group $W$ acts transitively on the
coroots [H, $\S$10.4, Lemma C].  Thus,
$\mu (\alpha_i^\vee) \in 2 \Z \cap 3 \Z = 6 \Z$ for all the simple coroots
$\alpha_i^\vee$.  This proves that $\mu \in 6\Lambda$, i.e.,
$$[M_6]_W \subset 6\Lambda . \leqno{(2)} $$
Comparing (1) and (2) we get $[M_6]_W=6\Lambda$.
This proves the theorem in the case of $E_6$.

\paragraph*{h) $E_7$:} By the chart and the result for the cases of $A_\ell,
D_\ell$, we get
$$
\begin{array}{lll}
L(E_7; \alpha_1, \dots, \alpha_7) &=& [(\Z \theta+L( D_6; \alpha_7, \alpha_6,
\dots, \alpha_2)) \cap \\
&&L(A_7; - \theta, \alpha_1, \alpha_3, \alpha_4, \dots, \alpha_7) \cap \\
&&(L(A_2; -\theta, \alpha_1)+ L(A_5; \alpha_2, \alpha_4, \alpha_5,
\alpha_6, \alpha_7))\cap \\
 &&(L(A_5; -\theta, \alpha_1, \alpha_3, \alpha_4, \alpha_2)+L(A_2; 
\alpha_6,\alpha_7))
\cap\\
&&(L(D_6; - \theta, \alpha_1, \alpha_3, \alpha_4, \alpha_5, \alpha_2)+ \Z
\alpha_\ell)]_W \\
&=&[(\Z \theta + L' (\alpha_7, \alpha_6, \dots , \alpha_2)) \cap \\
&&(\Z \theta + \Z \alpha_1+ \Z \alpha_3+ \Z \alpha_4+ \cdots + \Z 
\alpha_7) \cap \\
&&(\Z \theta + \Z  \alpha_1+ \Z \alpha_2+ \Z \alpha_4 + \cdots + \Z \alpha_7)
\cap \\
&&(\Z \theta + \Z  \alpha_1 + \cdots + \Z \alpha_4+\Z \alpha_6+ \Z \alpha_7)
\cap \\
&&(L' (- \theta, \alpha_1, \alpha_3, \alpha_4, \alpha_5, \alpha_2)+ \Z
\alpha_7)]_W \\
&=& [M_7]_W, 
\end{array}
$$
where
$M_7 := \Z4 \alpha_1+ \Z 2 \alpha_2+ \Z 6 \alpha_3+ \Z 2 \alpha_4+ \Z 6 \alpha_5+ \Z 4 
\alpha_6+ \Z 2 \alpha_7.$
 
Clearly,
$12 \Lambda \subset M_7$ and since $\Lambda$ is $W$-stable, $12 \Lambda \subset [M_7]_W$.  

For any $\mu \in [M_7]_W$, by considering $(w \mu) (\alpha_1^\vee)$ and $
(w \mu) (\alpha_3^\vee)$ as in the proof of the theorem for the case of $E_6$, we get 
that
$\mu (\alpha_i^\vee ) \in 4 \Z \cap 6 \Z = 12 \Z$ for all the simple
coroots $\alpha_i^\vee$.
Thus,
$$ [M_7]_W \subset 12 \Lambda \ {\rm and \ hence} \ [M_7]_W=12\Lambda.$$
This takes care of the case of $E_7$.

Finally, we come to the following:
\paragraph*{(i) ~$ E_8$:} By the chart and the theorem for $E_6, E_7, D_8$ and
$A_\ell$, we get (denoting the weight lattice of  $E_i$ by  $\Lambda(E_i))$.
$$
\begin{array}{lll}
L(E_8; \alpha_1, \dots , \alpha_8) &=& [L(D_8; - \theta, \alpha_8, \alpha_7,
\dots, \alpha_2) \cap \\
&&L (A_8; \alpha_1, \alpha_3, \alpha_4, \dots,
\alpha_8, - \theta) \cap \\
&&(L (A_4; \alpha_1, \alpha_3, \alpha_4, \alpha_2) +L (A_4; \alpha_6, \alpha_7,
\alpha_8, - \theta)) \cap \\
&&(L(E_6; \alpha_1, \dots, \alpha_6)+ L (A_2; \alpha_8, - \theta)) \cap \\
&&(L(E_7; \alpha_1, \dots, \alpha_7)+\Z \theta)]_W\\
&=& [L' (\theta, \alpha_8, \alpha_7, \dots, \alpha_2) \cap \\
&&(\Z \alpha_1+ \Z \alpha_3+\Z \alpha_4+ \cdots + \Z \alpha_8+ \Z \theta) \cap \\
&&(\Z \alpha_1+ \cdots + \Z \alpha_4 + \Z \alpha_6+\Z \alpha_7+ \Z \alpha_8+\Z
\theta) \cap \\
&&(6 \Lambda(E_6)+ \Z \alpha_8+\Z \theta) \cap 
(12 \Lambda (E_7)+\Z \theta)]_W \\
&=& [M'_8 \cap (6 \Lambda(E_6)+ \Z \alpha_8+\Z \theta) \cap (12 \Lambda (E_7)+\Z 
\theta)]_W, 
\end{array}
$$
where
$$M'_8:=\{\Z 4 \alpha_1+ n3 \alpha_2+m \alpha_3+\Z 2 \alpha_4+\Z 10 \alpha_5+\Z
2\alpha_6+ \Z 2\alpha_7 
+ \Z 2 \alpha_8: n, m \in \Z \ {\rm and} \ n+m \ {\rm is \ even}\}.$$

But the coefficient of $\alpha_3$ in any element of $6\Lambda(E_6)+\Z \alpha_8+\Z
\theta$ is even. Thus,
$$[M'_8 \cap  (6\Lambda (E_6)+ \Z \alpha_8+ \Z \theta) \cap (12 \Lambda(E_7)+\Z \theta)]_W=$$
$$[M_8 \cap (6 \Lambda (E_6)+ \Z \alpha_8+\Z \theta) \cap (12 \Lambda (E_7)+ \Z \theta)]_W,$$
where
$$M_8:= \Z 4 \alpha_1+\Z 6 \alpha_2+\Z 2 \alpha_3+ \Z 2 \alpha_4+\Z 10 \alpha_5+
\Z 2\alpha_6+ \Z 2\alpha_7+ \Z 2 \alpha_8.$$
For any $\mu \in [M_8]_W$, by considering
$(w \mu) (\alpha_1^\vee), (w \mu)(\alpha_2^\vee)$ and $(w \mu) (\alpha_5^\vee)$ as in the 
proof of the theorem for the case of $E_6$, we get
$$[M_8]_W \subset 60 \Lambda (E_8). \leqno{(3)}$$
Conversely, since $\Lambda (E_8) = Q (E_8)$, we get that
$$60\Lambda (E_8) \subset M_8 \cap (6\Lambda (E_6)+\Z \alpha_8+ \Z \theta) \cap
(12 \Lambda (E_7)+\Z \theta),
\leqno{(4)}$$
 and hence
$$60 \Lambda(E_8) \subset [M_8 \cap (6\Lambda (E_6)+\Z \alpha_8+\Z \theta) \cap
(12 \Lambda (E_7)+\Z \theta)]_W.$$
To prove (4), it suffices to show that
$$ 60 \alpha_7 \in 6 \Lambda (E_6)+\Z \alpha_8+\Z \theta \leqno{(5)}$$
and
$$ 60 \alpha_8 \in 12 \Lambda (E_7)+\Z \theta. \leqno{(6)}$$
To prove (5), observe that
\begin{align*} 60 \alpha_7 &=20 \theta - 40 \alpha_8 - 20 (2 \alpha_1+ 3 \alpha_2+4 
\alpha_3
 +6 \alpha_4 + 5 \alpha_5+ 4 \alpha_6)\\
&=20  \theta - 40 \alpha_8 - 60 \omega_6 (E_6),
\end{align*}
where $\omega_6 (E_6)$ is the sixth fundamental weight of $E_6$.  Similarly, to
prove (6), observe that
\begin{align*} 60 \alpha_8 &=30 \theta - 30 (2 \alpha_1+ 3\alpha_2+ 4 \alpha_3+ 6 
\alpha_4+ 5 \alpha_5+ 4 \alpha_6+ 3 \alpha_7)\\
&= 30 \theta - 60 \omega_7 (E_7).
\end{align*}
This proves (5) and (6) and thus (4).  Combining (3) and (4), we get that
$$[M_8 \cap (6 \Lambda (E_6)+\Z \alpha_8+\Z \theta) \cap (12 \Lambda(E_7)+ \Z \theta)]_W =
60 \Lambda (E_8)=60Q.$$
This proves the theorem for $E_8$ and hence the theorem is completely proved.\\
$~~~~~~~~~~~\hfill{\Box}$

\begin{remark} Theorem 3.10 (a) is obtained earlier by B. Howard [Ho].
\end{remark}

\vskip5ex

\noindent
Address: Department of Mathematics, UNC at Chapel Hill, Chapel Hill, 
NC 27599-3250, USA (email:  shrawan@email.unc.edu)
\end{document}